\documentclass[12pt,reqno]{amsart}
\usepackage{amssymb,amsmath,graphicx,
	amsfonts,euscript}
\usepackage{color}
\usepackage{mathrsfs}
\theoremstyle{plain}

\numberwithin{equation}{section}

\begin{document}

\title[Global bifurcation structure and geometric properties]{ Global bifurcation structure and geometric properties for steady periodic water waves with vorticity}

\author[Guowei Dai, Yong Zhang]{Guowei Dai, Yong Zhang}

\address[Guowei Dai]{School of Mathematical Sciences, Dalian University
of Technology, Dalian, 116024, China}

\address[Yong Zhang]{ School of Mathematical Sciences, Dalian University
of Technology, Dalian, 116024, China} \email{18842629891@163.com}

\begin{abstract}
This paper studies the classical water wave problem with vorticity described by the Euler equations
with a free surface under the influence of gravity over a flat bottom. Based on fundamental work
\cite{ConstantinStrauss}, we first obtain two continuous bifurcation curves which meet the laminar
flow only one time by using modified analytic bifurcation theorem. They are
symmetric waves whose profiles are monotone between each crest and trough. Furthermore, we find
that there is at least one inflection point on the wave profile between successive crests and troughs and the free surface is strictly concave at any crest and strictly
convex at any trough. In addition, for favorable vorticity, we prove that the vertical
displacement of water waves decreases with depth.
	\end{abstract}
\maketitle

%\noindent {\bf Mathematics Subject Classification (2010)} {76A10; 76D03 }
\section{Introduction}

\bigskip
\quad\,

It's often possible for us to observe water wave while watching the sea or a lake. In fact, the systematic study of irrotational water waves can date back to the 1840s. In 1847, Stokes \cite{Stokes1} conjectured the existence of a large amplitude periodic wave with a stagnation point and a corner containing an angle of 120$^{\circ}$ at its highest point.  The first rigorous constructions by power
series of such waves due to Nekrasov \cite{Nekrasov}, which were local in the sense that
the wave profiles were almost flat. Constructions of large-amplitude irrotational
waves were begun by Krasovskii \cite{Krasovskii} which were refined by Keady and Norbury \cite{Keady} by using the methods of global bifurcation theory.
Buffoni, Dancer and Toland \cite{BDT1, BDT2}, Buffoni and Toland \cite{BT} used the global analytic bifurcation theory
to obtain the existence of waves of all amplitudes from zero up to that of Stokes' highest wave. Besides, Amick \cite{Amick} also proved that for any irrotational wave the angle $\theta$ of inclination (with respect to the horizontal) of the profile must be less than 31.15$^{\circ}$.

In recent years, the classical hydrodynamic problem concerning two-dimensional steady periodic travelling water waves with vorticity has attracted considerable interests, starting with the study of Constantin and Strauss\cite{ConstantinStrauss} for periodic waves of finite depth. In 2004, assuming there is no stagnation points, Constantin and Strauss\cite{ConstantinStrauss} obtained two-dimensional inviscid periodic traveling waves with vorticity by using bifurcation theory.
By local bifurcation theory, Constantin and V$\breve{a}$rv$\breve{a}$ruc$\breve{a}$ \cite{Constantin1} studied periodic traveling gravity waves with stagnation points at the free surface of water in a flow of constant vorticity over a flat bed. Later, Constantin, Strauss and V$\breve{a}$rv$\breve{a}$ruc$\breve{a}$ \cite{Constantin0} further extended the result to obtain global bifurcation curve with critical
layers via analytic global bifurcation theorem which was obtained by Dancer \cite{D1} and improved by Buffoni
and Toland \cite[Theorem 9.1.1]{BT}. Recently, Constantin, Strauss and V$\breve{a}$rv$\breve{a}$ruc$\breve{a}$ \cite{Constantin2} strictly proved that the downstream waves on a global bifurcation branch are never overhanging and show that the waves approach their maximum possible amplitude. The authors \cite{Dai} investigated the local behavior of the bifurcation curve at the bifurcation point, determined the linearized stability of the curve when the vorticity is small, and obtained some further qualitative information about the global behavior of the curve. On the other hand, Strauss et al. in \cite{Strauss1,Strauss2} established an upper bound of 45$^{\circ}$ on $\theta$ for a large class of waves with favorable and small adverse vorticity.

Note that Constantin and Strauss \cite{ConstantinStrauss} obtain a global connected set of water waves of large amplitude without stagnation points. However, the global structure and some geometric properties of this set are not very clear. The main aim of this paper is to study global structure and some geometric properties of the connected set obtained in \cite{ConstantinStrauss}.  In Section 2, we recall the governing equations for two-dimensional steady  periodic water waves in different formulations and state our main results. The Section 3 is devoted to obtaining two global bifurcation solution curves and analysing its  structure. In Section 4, we mainly obtain some geometric properties of water waves to finish the proof of main theorem.

\section{Preliminaries and main results}

We recall the governing equations for the propagation of two-dimensional
gravity water waves. Choosing coordinates $(x, y)\in \mathbb{R}^2$ so that the horizontal $x$-axis is in the direction
of wave propagation, the $y$-axis points vertically upwards, and the origin lies at the
mean water level. In its undisturbed state the equation of the flat surface
is $y = 0$, and the flat bottom is given by $y =-d$ for some $d > 0$. In the presence
of waves, let $y = \eta(t, x)$ be the free surface and let $(u(t, x, y), v(t, x, y))$ be the
velocity field. Then the problem of traveling gravity water waves in $\mathbb{R}^2$ can be formulated as
\begin{equation}\label{waterwavetime}
\left\{
\begin{array}{ll}
u_x+v_y=0,\\
u_t+uu_x+vu_y=-P_x,\\
v_t+uv_x+vy_y=-P_y-g, \\
P=P_{\text{atm}}\,\,\,&\text{on}\,\,y=\eta(t,x),\\
v=\eta_t+u\eta_x &\text{on}\,\,y=\eta(t,x),\\
v=0 &\text{on}\,\, y=-d,
\end{array}
\right.
\end{equation}
where $P(t, x, y)$ denotes the pressure, $g$ is the gravitational constant of acceleration and $P_{\text{atm}}$ being the constant atmospheric pressure.

Given $c > 0$, we are looking for periodic waves traveling at speed $c$. So, the space-time dependence of the free surface, of the pressure, and of the velocity field has the form $(x -ct)$.
Define the stream function $\psi(x, y)$ by $\psi_x=-v$, $\psi_y=u-c$.
Then $-\Delta\psi=\gamma(\psi)$ where $\gamma$ is the vorticity function.
The relative mass flux is defined by
\begin{equation}
\int_{-d}^{\eta(x)}\psi_y\,dy:=p_0,\nonumber
\end{equation}
which is a negative constant if $u<c$.
From now on we always assume there is no stagnation points throughout the
fluid, that is to say $u<c$.
Let
\begin{equation}
\Gamma(p)=\int_0^p\gamma(-s)\,ds\nonumber
\end{equation}
have minimum value $\Gamma_0$ for $p_0\leq p\leq 0$.

Define $\overline{\Omega_\eta}$ as the closure of the open fluid domain
\begin{equation}
\Omega_\eta=\left\{(x,y)\in \mathbb{R}^2:x\in \mathbb{R},-d<y<\eta(x)\right\} \nonumber
\end{equation}
and let
\begin{equation}
\Omega_+=\left\{(x,y)\in \mathbb{R}^2:x\in (0,\pi),-d<y<\eta(x)\right\} \nonumber
\end{equation}
For $m \geq1$ an integer and $\alpha\in (0, 1)$, a domain $\Omega\subseteq  \mathbb{R}^2$ is a $C^{m+\alpha}$ domain
if each point of its boundary $\partial \Omega$ has a neighborhood in which $\partial \Omega$ is the graph
of a function with H\"{o}lder-continuous derivatives (of exponent $\alpha$) up to order $m$. Given a $C^{m+\alpha}$
domain $\Omega\subseteq  \mathbb{R}^2$, define
the space $C^{m+\alpha}_{2\pi} \left(\overline{\Omega}\right)$ of functions $f : \overline{\Omega} \rightarrow \mathbb{R}$ with H\"{o}lder-continuous derivatives
(of exponent $\alpha$) up to order $m$ and satisfying an $2\pi$-periodicity condition in the
$x$-variable. We use a similar notation for the case $\alpha= 0$ and for H\"{o}lder spaces of
functions of one variable.

 From Bernoulli's law, we know that
\begin{equation}
E=\frac{(c-u)^{2}+v^{2}}{2}+gy+P+\Gamma(-\psi) \nonumber
\end{equation}
is a constant along each streamline. Therefore, the dynamic boundary condition is equivalent to
\begin{equation}
|\nabla\psi|^{2}+2g(y+d)=Q~~~~~~on ~~y=\eta(x),  \nonumber
\end{equation}
where $Q=2(E-P_{atm}+gd)$ is a constant. According to \cite{Constantin2011} and \cite{ConstantinStrauss}, then the problem (\ref{waterwavetime}) can be reformulated as
\begin{equation}\label{waterwave}
\left\{
\begin{array}{ll}
\Delta \psi=-\gamma(\psi)\,\, &\text{in}\,\, -d<y<\eta(x),\\
\vert \nabla \psi\vert^2+2g(y+d)=Q & \text{on}\,\,y=\eta(x),\\
\psi=0 &\text{on}\,\,y=\eta(x),\\
\psi=-p_0 &\text{on}\,\, y=-d,
\end{array}
\right.
\end{equation}
where the constant $Q/2g$ is the total head. Under the assumption $u<c$, let $q=x$, $p=-\psi$ and $h(q,p)=y+d$ (see \cite{ConstantinStrauss}), then problem (\ref{waterwave}) is equivalent to
\begin{equation}\label{hfunction}
\left\{
\begin{array}{ll}
\left(1+h_q^2\right)h_{pp}-2h_ph_qh_{pq}+h_p^2h_{qq}=-\gamma(-p)h_p^3\,\, &\text{in}\,\, p_0<p<0,\\
1+h_q^2+(2gh-Q)h_p^2=0&\text{on}\,\, p=0,\\
h=0 &\text{on}\,\, p=p_0,
\end{array}
\right.
\end{equation}
with $h$ even and of period $2\pi$ in the $q$-variable.

From \cite[Lemma 3.2]{ConstantinStrauss} we know that
problem (\ref{hfunction}) has the trivial solution
\begin{equation}
H(p)=\int_0^p\frac{1}{\sqrt{\lambda+2\Gamma(s)}}\,ds\nonumber
\end{equation}
for $\lambda\in\left[-2\Gamma_0,Q\right)$ with
\begin{equation}
Q=\lambda+2g\int_{p_0}^0\frac{1}{\sqrt{\lambda+2\Gamma(s)}}\,ds.\nonumber
\end{equation}
And the linearized problem of (\ref{hfunction}) at $H$ is
\begin{equation}\label{linearizedproblem}
\left\{
\begin{array}{ll}
m_{pp}+H_p^2m_{qq}=-3\gamma(-p)H_p^2m_p\,\, &\text{in}\,\, p_0<p<0,\\
gm=\lambda^{\frac{3}{2}}m_p&\text{on}\,\, p=0,\\
m=0 &\text{on}\,\, p=p_0,
\end{array}
\right.
\end{equation}
with $m$ even and $2\pi$-periodic in $q$.
From \cite[Lemma 3.3]{ConstantinStrauss}, there exists $\lambda^*>-2\Gamma_0$ and a solution $m(q,p)\not\equiv0$ of
(\ref{linearizedproblem}) that is even and $2\pi$-periodic in $q$.

Let $R$ be the rectangle $(-\pi,\pi)\times \left(p_0,0\right)$, $T=\{p=0\}$ be the top and $B=\left\{p=p_0\right\}$ the bottom of its closure $\overline{R}$.
Let
\begin{equation}
X=\left\{h\in C_{\text{per}}^{3+\alpha}\left(\overline{R}\right):h=0\,\,\text{on}\,\,B\right\},\,\,Y=Y_1\times Y_2=C_{\text{per}}^{1+\alpha}\left(\overline{R}\right)\times C_{\text{per}}^{2+\alpha}\left(\overline{R}\right),\nonumber
\end{equation}
where the subscript ''per'' means $2\pi$-periodic and even in $q$.
For any $\delta>0$, set
\begin{equation}
\mathcal{O}_\delta=\left\{(Q,h)\in \mathbb{R}\times X:h_p>\delta\,\, \text{in}\,\, \overline{R}, h<\frac{Q-\delta}{2g}\,\,\text{on}\,\, T\right\}.\nonumber
\end{equation}
Define
\begin{equation}
F=\left(F_1,F_1\right):\mathbb{R}\times X\rightarrow Y\nonumber
\end{equation}
with
\begin{equation}
F_1(Q,h)=\left(1+h_q^2\right)h_{pp}-2h_ph_qh_{pq}+h_p^2h_{qq}+\gamma(-p)h_p^3\nonumber
\end{equation}
and
\begin{equation}
F_2(Q,h)=1+h_q^2+(2gh-Q)h_p^2.\nonumber
\end{equation}
Thus we have that
\begin{equation}
F(Q(\lambda),H(p))=0,\nonumber
\end{equation}
where
\begin{equation}
Q(\lambda)=\lambda+2g\int_{p_0}^0\frac{1}{\sqrt{\lambda+2\Gamma(s)}}\,ds.\nonumber
\end{equation}
Constantin and Strass \cite{ConstantinStrauss} established the following milestone result.
\\ \\
\textbf{Theorem A.} \emph{Let the wave speed $c > 0$, the wavelength $2\pi$,
and the relative mass flux $p_0 < 0$ be given. For a constant $\alpha\in (0, 1)$, let the
vorticity function $\gamma\in C^{1+\alpha}\left[0, \left\vert p_0\right\vert\right]$ be such that}
\begin{equation}
\int_{p_0}^0\left[\left(p-p_0\right)^2\left(2\Gamma(p)-2\Gamma_0\right)^{\frac{1}{2}}
+\left(2\Gamma(p)-2\Gamma_0\right)^{\frac{3}{2}}\right]\,dp<gp_0^2.\nonumber
\end{equation}
\emph{Consider traveling solutions of speed c and relative mass flux $p_{0}$ of the water
wave problem (\ref{waterwavetime}) with vorticity function $\gamma$ such that $u < c$ throughout the
fluid. There exists a connected set $\mathcal{C}$ of solutions $(u, v, \eta)$ in the space $C^{2+\alpha}_{per} \left(\overline{\Omega_\eta}\right)\times
C^{2+\alpha}_{per} \left(\overline{\Omega_\eta}\right)\times C^{3+\alpha}_{per} (\mathbb{R})$ with the following proper}
\begin{itemize}
  \item \emph{$\mathcal{C}$ contains a trivial laminar flow};
  \item \emph{there is a sequence $\left(u_n, v_n, \eta_n\right)\in \mathcal{C}$ for which $\max_{\overline{\Omega_{\eta_n}}}u_n\uparrow c$}.
\end{itemize}
\emph{Furthermore, each nontrivial solution $(u, v, \eta) \in \mathcal{C}$ satisfies the following}:

(i) \emph{$u$, $v$ and $\eta$ have period $2\pi$ in $x$};

(ii) \emph{within each period the wave profile $\eta$ has a single maximum and a
single minimum; say the maximum occurs at $x = 0$};

(iii) \emph{$u$ and $\eta$ are symmetric while $v$ is antisymmetric around the line $x = 0$};

(iv) \emph{a water particle located at $(x, y)$ with $0 < x < \pi$ and $y > -d$ has positive
vertical velocity $v > 0$};

(v) \emph{$\eta'(x) < 0$ on $(0, \pi)$}.
\\
\\
Here we further improve Theorem A by showing the following theorem.
\\ \\
\textbf{Theorem 2.1.} \emph{Under the assumptions of Theorem A, there exist two continuous curves $\mathcal{K}^+$ and $\mathcal{K}^-$ of solutions $(u, v, \eta)$ in the space $C^{2+\alpha}_{per} \left(\overline{\Omega_\eta}\right)\times
C^{2+\alpha}_{per} \left(\overline{\Omega_\eta}\right)\times C^{3+\alpha}_{per} (\mathbb{R})$ with the following properties}
\begin{itemize}
  \item \emph{$\mathcal{K}^\nu$ ($\nu\in\{+,-\}$) contains the trivial laminar flow if and only if $Q=Q^*$ for some positive $Q^*>0$};
  \item \emph{there is a sequence $\left(u_n, v_n, \eta_n\right)\in \mathcal{K}^\nu$ for which $\max_{\overline{\Omega_{\eta_n}}}u_n\uparrow c$}.
\end{itemize}
\emph{Furthermore, each nontrivial solution $(u, v, \eta) \in \mathcal{K}^\nu$ satisfies the properties (i)-(v) described in Theorem A. In addition, there also holds}:

(vi) \emph{there is at least one inflection point on the free surface $y=\eta(x)$ for $x\in[0,\pi]$, and the free surface is strictly concave at any crest and strictly convex at any trough (see Figure 1)};

\begin{figure}[ht]
\centering
\includegraphics{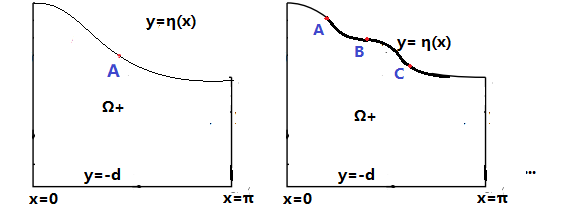}
\caption{The free surface with odd inflection points}
\label{fig3}
\end{figure}

(vii) \emph{For favorable vorticity (i.e $\gamma\leq 0, \gamma'\leq 0$), the vertical displacement $L(y)$ of water waves on any streamline decreases with depth (see Figure 2)}.

\begin{figure}[ht]
\centering
\includegraphics{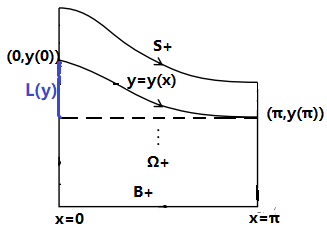}
\caption{The vertical displacement $L(y)$ of water waves on the streamline}
\label{fig4}
\end{figure}

Compared with Theorem A, we obtain two continuous curves which meet the laminar flow only one time.
Note that A. Constantin \cite[Theorem 3.5]{Constantin2011} also obtained a continuous curve $\mathscr{K}_\delta$, which may form a closed loop.
Here $\mathcal{K}^\nu$ does not form a closed loop. Besides, we also obtain new geometric properties $(vi)$ and $(vii)$ of the solution sets.

\section{The global bifurcation structure}
\bigskip
\quad\,
To prove Theorem 2.1, we first improve the analytic global bifurcation theorem of Buffoni
and Toland \cite[Theorem 9.1.1]{BT} as follows.
\\ \\
\textbf{Theorem 3.1.} \emph{Let $X$ and $Y$ be Banach spaces, $\mathcal{O}$ be
an open subset of $\mathbb{R}\times X$ and $F: \mathcal{O}\rightarrow Y$ be a real-analytic function. Suppose that}

(H1) \emph{$F(\lambda, 0)=0$ for all $(\lambda, 0)\in \mathcal{O}$};

(H2) \emph{for some $\lambda_*\in \mathbb{R}$, $\mathcal{N} \left(\partial_u F\left(\lambda_*, 0\right)\right)$ and $Y /\mathcal{R}\left(\partial_u F\left(\lambda_*, 0\right)\right)$ are $1$-dimensional, with
the null space generated by $u_*$, and the transversality condition}
\begin{equation}
\partial_{\lambda,u}^2F\left(\lambda_*, 0\right)\left(1,u_*\right)\not\in\mathcal{R}\left(\partial_u F\left(\lambda_*, 0\right)\right)\nonumber
\end{equation}
\emph{holds, where $\mathcal{N} \left(\partial_u F\left(\lambda_*, 0\right)\right)$ and $\mathcal{R}\left(\partial_u F\left(\lambda_*, 0\right)\right)$ denote null space and range space of $\partial_u F\left(\lambda_*, 0\right)$, respectively};

(H3) \emph{$\partial_uF(\lambda, u)$ is a Fredholm operator of index zero for any $(\lambda, u)\in \mathbb{O}$ such that
$F(\lambda, u)=0$};

(H4) \emph{all bounded closed subsets of $\left\{(\lambda,u)\in \mathcal{O}:F(\lambda,u)=0\right\}$ are compact in $\mathbb{R}\times X$}.

\emph{Then there exist in $\mathcal{O}$ two continuous curve $\mathcal{K}^\nu=\left\{(\lambda(s), u(s)):\nu s\geq0\right\}$ ($\nu\in\{+,-\}$) of solutions to
$F(\lambda, u)=0$ such that}:

(C1) \emph{$(\lambda(0), u(0))=\left(\lambda_*, 0\right)$};

(C2) \emph{$u(s)=su_*+o(s)$ in $X$, $\vert s\vert<\varepsilon$ as $s\rightarrow 0$};

(C3) \emph{there exist a neighbourhood $\mathcal{W}$ of $\left(\lambda_*, 0\right)$ and $\varepsilon>0$ sufficiently small such that}
\begin{equation}
\left\{(\lambda,u)\in \mathcal{W}:u\neq0\,\,\text{and}\,\,F(\lambda,u)=0\right\}=\left\{(\lambda(s),u(s)):0<\vert s\vert<\varepsilon\right\};\nonumber
\end{equation}

(C4) \emph{$\mathcal{K}^\nu$ has a real-analytic reparametrization locally around each of its points};

(C5) \emph{one of the following alternatives occurs}:

(1) \emph{$(\lambda(s),u(s))\rightarrow\infty$ in $\mathbb{R}\times X$ as $s\rightarrow\infty$};

(2) \emph{$(\lambda(s),u(s))$ approaches $\partial \mathcal{O}$ as $s\rightarrow\infty$};

(3) \emph{$\mathcal{K}^\nu$ contains a trivial point $(\mu,0)\in \mathcal{O}$ with $\mu\neq\lambda_*$}.

\emph{Moreover, such a curve of solutions to $F(\lambda, u)=0$ having the properties (C1)--(C5)
is unique (up to reparametrization)}.
\\

In \cite[Theorem 9.1.1]{BT} or \cite[Theorem 6]{Constantin0}, $\mathcal{K}^\nu$ may form a closed loop, that is, there exists $T >0$ such that $(\lambda(s+T), u(s+T))=(\lambda(s), u(s))$ for all $s\in \mathbb{R}$.
Here we further show that $\mathcal{K}^\nu$ contains another bifurcation point if a closed loop occurs, which is  coincident with the Rabinowitz Global Bifurcation Theorem \cite{Rabinowitz0}.
\\ \\
\textbf{Proof of Theorem 3.1.} A distinguished arc is a maximal connected subset of $\mathcal{K}^+$. Let $\mathcal{A}_0$ denote the first distinguished arc of $\mathcal{K}^+$ which bifurcates from $\left(\lambda_0,0\right)$. Here we only show the case of $\nu=+$ for simplicity. Suppose by contradiction that
$\mathcal{K}^+$ forms a closed loop and does not contain another bifurcation point. From \cite[Theorem 9.1.1]{BT} there exists some constant $T>0$ such that
$(\lambda(T),u(T))=\left(\lambda_0,0\right)$. So there is a segment of $\mathcal{K}^+$, parameterized by $s<T$ sufficiently
close to $T$, which is a subset of $\left\{(\lambda(s),u(s)):0<s<\varepsilon\right\}:=\mathcal{S}^+$.
Since $\mathcal{S}^+\subseteq \mathcal{A}_0$, there exist positive sequences $\left\{s_k\right\}$ and $\left\{t_k\right\}$ with $s_k\downarrow 0$ and $t_k\downarrow0$ such that
\begin{equation}
\left(\lambda\left(s_k\right),u\left(s_k\right)\right)=\left(\lambda\left(T-t_k\right),u\left(T-t_k\right)\right).\nonumber
\end{equation}
By \cite[Theorem 9.1.1]{BT} $T-s_k-t_k$ is an integer multiple of $T$.
That is to say there exists a $m\in \mathbb{N}$ such that
$(1-m)T=s_k+t_k$ for all $k$. It follows that $m=1$. Thus, $s_k+t_k=0$ for all $k$, which is a contradiction.\qed
\\

\indent Applying Theorem 3.1 to $F(Q,h)$ on $\mathcal{O}_\delta$, we obtain the following theorem.
\\ \\
\textbf{Theorem 3.2.} \emph{There exist two continuous curve $\mathcal{K}_\delta^\nu=\left\{(Q(s), h(s)):\nu s\geq0\right\}$ ($\nu\in\{+,-\}$) of solutions to $F(Q, h)=0$. Either}

(i) \emph{$\mathcal{K}_\delta^\nu$ is unbounded in $\mathbb{R}\times X$, or}

(ii) \emph{$\mathcal{K}_\delta^\nu$ contains a point $(Q,h)\in \partial\mathcal{O}_\delta$}.

\noindent \emph{Moreover, $\mathcal{K}_\delta^\nu\cap \left\{(Q,H(p))\right\}=\left(Q^*,H^*\right)$ where
$Q^*=Q\left(\lambda^*\right)$ and $H^*=H\left(\lambda^*\right)$}.
\\ \\
\textbf{Proof of Theorem 3.2.} From \cite[Lemma 3.7, 3.8 and Theorem 3.1]{ConstantinStrauss} we have that $F$ satisfies (H2) in Theorem 3.1.
By \cite[Lemma 4.1 and 4.3]{ConstantinStrauss}, $F$ satisfies (H2)-(H3) in Theorem 3.1.

So, by Theorem 3.1, there exist two continuous curve $\mathcal{K}_\delta^\nu=\left\{(Q(s), h(s)):\nu s\geq0\right\}$ ($\nu\in\{+,-\}$) of solutions to $F(Q, h)=0$ such that

(C1) $(Q(0), h(0))=\left(Q\left(\lambda^*\right), H(p)\right)$;

(C2) $h(s)=H(p)+s\varphi_*+o(s)$ in $X$, $\vert s\vert<\varepsilon$ as $s\rightarrow 0$, where $\varphi_*\in \mathcal{N} \left(\partial_h F\left(Q\left(\lambda^*\right), H(p)\right)\right)$ with $\left \Vert \varphi_*\right\Vert_X=1$;

(C3) there exist a neighbourhood $\mathcal{W}$ of $\left(Q\left(\lambda^*\right), H(p)\right)$ and $\varepsilon>0$ sufficiently small such that
\begin{equation}
\left\{(Q,h)\in \mathcal{W}:h\neq H(p)\,\,\text{and}\,\,F(Q,h)=0\right\}=\left\{(Q(s),h(s)):0<\vert s\vert<\varepsilon\right\};\nonumber
\end{equation}

(C4) $\mathcal{K}_\delta^\nu$ has a real-analytic reparametrization locally around each of its points;

(C5) one of the following alternatives occurs:

(1) $(Q(s),h(s))\rightarrow\infty$ in $\mathbb{R}\times X$ as $s\rightarrow\infty$;

(2) $(Q(s),h(s))$ approaches $\partial \mathcal{O}_\delta$ as $s\rightarrow\infty$;

(3) $\mathcal{K}_\delta^\nu$ contains a trivial point $(Q(\mu),H(\mu))\in \mathcal{O}_\delta$ with $\mu\neq\lambda^*$.

\noindent By \cite[Lemmatas 5.1--5.3]{ConstantinStrauss} the alternative (3) is impossible.
So $\mathcal{K}_\delta^\nu$ is unbounded or reach to the boundary of $\mathcal{O}_\delta$.
\qed\\

The conclusions of Theorem 3.2 are better than the corresponding ones of \cite[Theorem 5.4]{ConstantinStrauss} or \cite[Theorem 3.5]{Constantin2011}.
Moreover, the argument here is more concise. Let $\mathcal{C}^\nu=\cup_{\delta>0}\mathcal{K}_\delta^\nu$.
From \cite[Section 7]{ConstantinStrauss}, there exists a sequence
$\left(Q_n,h_n\right)\in \mathcal{C}^\nu$ such that $\sup_{\overline{D_{\eta_n}}}u_n\rightarrow c$, where
$\left(u_n,v_n,\eta_n\right)$ is the solution of the water wave problem (\ref{waterwavetime}) corresponding to $h_n$.
Hence, from the definition of $h$ we see that $\left\Vert \left(h_n\right)_p\right\Vert_{C_{\text{per}}^{0}\left(\overline{R}\right)}\rightarrow+\infty$ as $n\rightarrow+\infty$.
It follows that $\left\Vert h_n\right\Vert_{C_{\text{per}}^{3+\alpha}\left(\overline{R}\right)}\rightarrow+\infty$ as $n\rightarrow+\infty$.
Therefore, $\mathcal{C}^\nu$ is unbounded in the direction of $X$.
\\ \\
\textbf{Nodal pattern of solutions on $\mathcal{K}^\nu$:} By the definition of $\mathcal{O}_\delta$, $\mathcal{K}_\delta^\nu$ increase as $\delta$ decreases.
So, $\mathcal{K}^\nu=\lim_{\delta\rightarrow0}\mathcal{K}_\delta^\nu$ and it is closed and continuous, which corresponds a solution set of problem (\ref{hfunction}).
It is enough to show that $\mathcal{K}^\nu\cap \left\{(Q,H(p))\right\}=\left(Q^*,H^*\right)$.
To do it, we study the nodal pattern of solutions on $\mathcal{K}^\nu$.
We only show the case of $\nu=+$ for simplicity.

If $D$ is the open rectangle $(0,\pi)\times\left(p_0,0\right)$,
we denote its sides by
\begin{equation}
\partial D_t=\{(q,0):q\in(0,\pi)\},\,\,\,\partial D_b=\left\{\left(q,p_0\right):q\in(0,\pi)\right\},\nonumber
\end{equation}
\begin{equation}
\partial D_l=\left\{(0,p):p\in\left(p_0,0\right)\right\},\,\,\, \partial D_r=\left\{\left(\pi,p\right):p\in\left(p_0,0\right)\right\}.\nonumber
\end{equation}
We consider the following properties of $h$
\begin{equation}\label{property1}
\left\{
\begin{array}{ll}
h_q<0\,\, &\text{in}\,\, D\cup \partial D_t,\\
h_{qp}<0&\text{on}\,\, \partial D_b,\\
h_{qq}<0 &\text{on}\,\, \partial D_l,\\
h_{qq}>0 &\text{on}\,\, \partial D_r,
\end{array}
\right.
\end{equation}
\begin{equation}\label{property2}
h_{qqp}\left(0,p_0\right)<0,\,\,\,h_{qqp}\left(\pi,p_0\right)>0.
\end{equation}
\begin{equation}\label{property3}
\text{either}\,\,h_{qq}(\pi,0)>0\,\,\text{or}\,\,h_{qqp}(\pi,0)<0; \text{either}\,\,h_{qq}(0,0)<0\,\,\text{or}\,\,h_{qqp}(0,0)>0.
\end{equation}
By \cite[Lemma 5.1]{ConstantinStrauss} properties (\ref{property1})--(\ref{property3}) hold in a small neighborhood of
$\left(Q^*,H^*\right)$ in $ \mathbb{R}\times X$ along the bifurcation curve $\mathcal{K}^+\setminus\left\{\left(Q^*,H^*\right)\right\}$.

Since $h_p>0$ on $\overline{R}$ and $h\in X$, we have that $\min_{\overline{R}}h_p:=k>0$.
So the condition (5.14) of Serrin's Maximum Principle \cite{ConstantinStrauss} is still valid for $h$ with $(Q,h)\in\mathcal{K}^+\setminus\left\{\left(Q^*,H^*\right)\right\}$.
Thus, as that of \cite[Lemma 5.2]{ConstantinStrauss}, the nodal properties (\ref{property1})--(\ref{property3})
hold along $\mathcal{K}^+\setminus\left\{\left(Q^*,H^*\right)\right\}$ unless there exists $\lambda\neq\lambda^*$
such that $(Q(\lambda),H(\lambda))\in \mathcal{K}^+$.

We \emph{claim} that $\mathcal{K}^+$ meets another bifurcation point is impossible.
If $(Q(\lambda),H(\lambda))\in \mathcal{K}^+$, there is a sequence of solutions
$\left(Q\left(\lambda_n\right),h_n\right)\in \mathcal{K}^+$ with $h_n\not\equiv H(p)$ such that
$\left(\lambda_n,h_n\right)\rightarrow\left(\lambda,H(\lambda)\right)$ in $\mathbb{R}\times X$.
Clearly, we have that
\begin{equation}
H_p=\frac{1}{\sqrt{\lambda+2\Gamma(p)}}\geq \frac{1}{\sqrt{\lambda+2\Gamma_1}}\nonumber
\end{equation}
for any $p\in\left[p_0,0\right]$, where $\Gamma_1=\max_{\left[p_0,0\right]} \Gamma(p)$.
It follows that
\begin{equation}
\left(h_n\right)_p\geq \frac{1}{2\sqrt{\lambda+2\Gamma_1}}\nonumber
\end{equation}
for $n$ large enough.
Differentiating $F_1$ we have that
\begin{eqnarray}
F_{1h}(Q,h)&=&\left(1+h_q^2\right)\partial_p^2-2h_ph_q\partial_p\partial_q+h_p^2\partial_q^2+2h_{pp}h_q\partial_q\nonumber\\
& &-2h_{pq}h_p\partial_q-2h_{pq}h_q\partial_p-2h_{pq}h_q\partial_p+2h_{qq}h_p\partial_p+2\gamma(-p)h_p^2\partial_p.\nonumber
\end{eqnarray}
So, the linear operators $F_{1h}\left(Q\left(\lambda_n\right),h_n\right)$ are uniformly elliptic (uniformly in $(q,p)$ and $n$) because
\begin{equation}
\left(1+\left(h_n\right)_q^2\right)\left(h_n\right)_p^2-\left(\left(h_n\right)_p\left(h_n\right)_q\right)^2
=\left(h_n\right)_p^2\geq\frac{1}{\left(\lambda+2\Gamma_1\right)}.\nonumber
\end{equation}
Differentiating $F_2$ we have that
\begin{eqnarray}
F_{2h}(Q,h)=2h_q\partial_q+2(2gh-Q)h_p\partial_p+2gh_p^2.\nonumber
\end{eqnarray}
From the definition of $H(p)$ with $p\in\left[p_0,0\right]$ we get that
\begin{equation}
\left\vert (Q-2gH)H_p\right\vert=\left(\lambda-2g\int_0^p\frac{1}{\sqrt{\lambda+2\Gamma(s)}}\,ds\right) \frac{1}{\sqrt{\lambda+2\Gamma_1}}\geq \frac{\lambda}{\sqrt{\lambda+2\Gamma_1}}.\nonumber
\end{equation}
It follows that
\begin{equation}
\left\vert \left(Q\left(\lambda_n\right)-2gh_n\right)\left(h_n\right)_p\right\vert\geq \frac{\lambda}{2\sqrt{\lambda+2\Gamma_1}}\nonumber
\nonumber
\end{equation}
for $n$ large enough.
Thus, the linear boundary operators $F_{2h}\left(Q\left(\lambda_n\right),h_n\right)$ are uniformly oblique.
Then as that of \cite[Lemma 5.3]{ConstantinStrauss} we deduce that $\lambda=\lambda^*$.\qed

\section{The geometric properties}

\bigskip
\quad\,
In this section, we aim to finish the proof of Theorem 2.1 by establishing the geometric properties $(vi)$ and $(vii)$.
\\

\textbf{The proof of Theorem 2.1.}
Suppose for seeking contradiction of $(vi)$ that there is no inflection point on $y=\eta(x)$ for $x\in[0,\pi]$. Then the free surface $y=\eta(x)$ must be strictly convex (or concave) for $x\in[0,\pi]$. Without loss of generality, we assume that it's strictly convex, i.e
\begin{equation}
\eta_{xx}>0~~~~for~~x\in(0,\pi). \label{eq4.1}
\end{equation}
Due to the symmetry and periodic properties of $v$ in Theorem A, we have that
\begin{equation}
v(\pi,y)=v(0,y)=0,~~~for~~y\in[-d, \eta(x)]. \label{eq4.2}
\end{equation}
Based on the mass conservation in (\ref{waterwavetime}), we obtain that
\begin{equation}
\frac{dv(x,\eta(x))}{dx}=v_{x}+v_{y}\eta_{x}=v_{x}-u_{x}\eta_{x}. \label{eq4.3}
\end{equation}
Besides, the fact $\psi(x,\eta(x))=0$ indicates that
\begin{equation}
\eta_{x}=-\frac{\psi_{x}}{\psi_{y}}=\frac{v}{u-c}<0, \label{eq4.4}
\end{equation}
the last inequality follows from property $(iv)$ in Theorem A.
Combining (\ref{eq4.1}) with (\ref{eq4.4}), we get that
\begin{equation}
\eta_{xx}=(\frac{v}{u-c})_{x}=\frac{v_{x}(u-c)-vu_{x}}{(u-c)^{2}}>0, \nonumber
\end{equation}
which means that
\begin{equation}
v_{x}(u-c)-vu_{x}>0. \label{eq4.5}
\end{equation}
From boundary condition in (\ref{waterwavetime}), (\ref{eq4.4}), (\ref{eq4.5}) and the property $(v)$ in Theorem A, we have
\begin{equation}
v\frac{dv(x,\eta(x))}{dx}=v(v_{x}+v_{y}\eta_{x})=(u-c)\eta_{x}(v_{x}-u_{x}\eta_{x})=\eta_{x}((u-c)v_{x}-vu_{x})<0, \nonumber
\end{equation}
which implies
\begin{equation}
\frac{dv(x,\eta(x))}{dx}<0,~~~~for ~~x\in(0,\pi). \nonumber
\end{equation}
This means $v$ decreases strictly along surface $y=\eta(x)$ from wave crest $(0,\eta(0))$ to trough $(\pi,\eta(\pi))$, which is contradicted with (\ref{eq4.2}). Indeed, the number of inflection points on surface $y=\eta(x)$ must be odd from wave crest $(0,\eta(0))$ to trough $(\pi,\eta(\pi))$ due to the fact (\ref{eq4.2}) and vertical velocity $v$ would decrease along the surface if the free surface is convex, otherwise $v$ would increase. For instance, if there are two inflection points on surface $y=\eta(x)$ from wave crest $(0,\eta(0))$ to trough $(\pi,\eta(\pi))$ (see Figure 3), then $0=v(O)<v(A), v(A)>v(B)>0$, $0<v(B)<v(C)=0$ is a contradiction.
\begin{figure}[ht]
\centering
\includegraphics{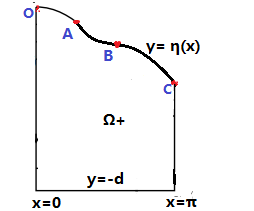}
\caption{The free surface with two inflection points}
\label{fig3}
\end{figure}

By the way, considering the stability of wave profile, we conjecture that there should be one inflection point on surface $y=\eta(x)$ from wave crest $(0,\eta(0))$ to trough $(\pi,\eta(\pi))$. In fact, this is also open problem for Stokes wave. As Constantin \cite[Section 4.4]{Constantin2011} mentioned that "The qualitative description of the flow beneath a smooth Stokes wave with no underlying current is almost complete. The only missing major aspect (other than understanding the wave of greatest height) is the behavior of the vertical velocity component $v$: we know how its sign depends on the location within the fluid domain, but how about its monotonicity? One would conjecture that along each streamline, between crest and trough, $v$ first increases with positive values away from the crest line and then decreases toward zero beneath the wave trough.", which means it was conjectured that for Stokes wave there should be one inflection point on surface $y=\eta(x)$ from wave crest $(0,\eta(0))$ to trough $(\pi,\eta(\pi))$. Here we only answer the question partially but our proof holds for any vorticity.

Now it remains to prove the property $(vii)$ in Theorem 2.1. Base on the result of Lemma 5.2 of Basu \cite{B}, we can state that for $\gamma\leq 0, \gamma'\leq 0$, the horizontal fluid velocity $u$ is a strictly decreasing function of $x$ along any streamline in $\overline{\Omega_{+}}$. (Note that the definition of vorticity function differs by a minus sign in this paper.) That is to say,
\begin{equation}
\frac{du(x,y(x))}{dx}=u_{x}+u_{y}y_{x}<0.  \label{eq4.6}
\end{equation}
In addition, $\psi(x,y(x))=c_{1}$ ($c_{1}$ is a constant) due to $y = y(x)$ is a streamline. We have $y_{x}=-\frac{\psi_{x}}{\psi_{y}}=\frac{v}{u-c} ~on ~y=y(x)$, hence
\begin{equation}
\frac{dy_{x}}{dy}=\frac{d(\frac{v}{u-c})}{dy}=\frac{v_{y}(u-c)-u_{y}v}{(u-c)^{2}}. \label{eq4.7}
\end{equation}
Thus, (\ref{eq4.6}), (\ref{eq4.7}) and the assumption $u<c$ imply that
\begin{equation}
\frac{dy_{x}}{dy}=\frac{-u_{x}(u-c)-u_{y}(u-c)y_{x}}{(u-c)^{2}}=-\frac{u_{x}+u_{y}y_{x}}{u-c}=\frac{u_{x}+u_{y}y_{x}}{c-u}<0. \label{eq4.8}
\end{equation}
Define the vertical displacement of water waves on streamline $y=y(x)$ by $L(y)$
\begin{equation}
L(y)=y(0)-y(\pi)=-\int^{\pi}_{0}y_{x}dx. \label{eq4.9}~
\end{equation}
By (\ref{eq4.8}) and (\ref{eq4.9}), we get
\begin{equation}
\frac{dL(y)}{dy}=\frac{d(y(0)-y(\pi))}{dy}=-\int^{\pi}_{0}\frac{dy_{x}}{dy}dx>0. \label{eq4.10}
\end{equation}
Up to now, we finish the proof of Theorem 2.1.
\qed
%Stokes conjectured
%
%\\ \\
%\textbf{Data Availability Statements.} Data sharing not applicable to this article as no datasets were generated or analysed during the current study.

%\\ \\
%\textbf{Acknowledgment}
%\bigskip\\
%\indent The first author thanks Academy of Mathematics and Systems Science of CAS for the
%invitation and hospitality during his visit.

\bibliographystyle{amsplain}
\makeatletter
\def\@biblabel#1{#1.~}
\makeatother

%\bibliography{mybib2}

\providecommand{\bysame}{\leavevmode\hbox to3em{\hrulefill}\thinspace}
\providecommand{\MR}{\relax\ifhmode\unskip\space\fi MR }
\providecommand{\MRhref}[2]{%
  \href{http://www.ams.org/mathscinet-getitem?mr=#1}{#2}
}
\providecommand{\href}[2]{#2}

\end{document}